\documentclass[11pt,leqno,oneside]{amsart}
\usepackage[width=6.5in,height=8.8in]{geometry}

\newtheorem{thm}{Theorem}[section]

\newtheorem{lem}[thm]{Lemma}
\newtheorem{cor}[thm]{Corollary}
\newtheorem{pro}[thm]{Proposition}
\newtheorem{ex}[thm]{Example}
\theoremstyle{definition}
\newtheorem{rem}[thm]{Remark}

\begin{document}
	\title{Lelong Numbers Above an Arbitrary  Bound and Plane Curves }
%\title{Upper level sets of Lelong numbers with smaller threshold}

\author{Al\.i Ula\c{s} \"Ozg\"ur K\.i\c{s}\.isel}
\address{Department of Mathematics, Middle East Technical University, 06800 Ankara, Turkey}
\email{akisisel@metu.edu.tr}

\author{Ozcan Yaz{\i}c{\i}}
\address{Department of Mathematics, Middle East Technical University, 06800 Ankara, Turkey} 
\email{oyazici@metu.edu.tr}
\date{\today}


\begin{abstract} Let $T$ be a positive closed current of bidimension $(1,1)$ with unit mass on $\mathbb P^2$ and $E^+_{\alpha}(T)$ be the upper level sets of Lelong numbers $\nu(T,x)$ of $T$. For any $\alpha\geq \frac{3}{d+5}$, $d\in \mathbb Z^+$ we show that  $|E^+_{\alpha}(T)\setminus C|\leq 2$ for a  curve $C$ of degree $d$. As a corollary of this result, we obtain a quantative lower bound for the integrability exponent at any point away from the curve $C$ and a two point set $P$.
\end{abstract}

\maketitle

\section{Introduction}
Suppose that $T$ is a positive closed current of bidimension $(1,1)$ on the complex projective plane $\mathbb{P}^2$ such that $\|T\|=1$. Let $\nu(T,p)$ denote the Lelong number of $T$ at a point $p\in \mathbb{P}^2$. For $0<\alpha<1$, set 
\begin{eqnarray*} 
E_{\alpha}(T):= \{p\in \mathbb{P}^{2} | \nu(T,p)\geq \alpha \}  \\ 
E^+_{\alpha}(T):= \{p\in \mathbb{P}^{2} | \nu(T,p) > \alpha \}.  \\ 
\end{eqnarray*}

Siu (\cite{Siu}) proved  that for any $\alpha >0$  $E_{\alpha}(T)$ is   an analytic subvariety of $\mathbb P^2$ of dimension at most $1$ . In \cite{Com}, Coman proved that if $\alpha \geq \frac{1}{2}$, $|E^+_{\alpha}(T)\setminus L|\leq 1$ for some complex line $L$  and if $\alpha \geq \frac{2}{5}$, $|E^+_{\alpha}(T)\setminus C|\leq 1$ for some conic. In \cite{KY} the authors showed that for any $\alpha \geq \frac{1}{3}$, $|E^+_{\alpha}(T)\setminus C|\leq 2$ for some cubic curve $C$. Some further geometric results are obtained for upper level sets in  \cite{Com, CT, Hef1, Hef2, KY2}.

In this paper our goal is to  investigate the geometry of this upper level set for arbitrarily smaller thresholds $\alpha$. 

For a finite subset $S$ of $\mathbb{P}^2$, the quantity $|S|$ will denote the number of points in $S$.
Our main result is the following:

\begin{thm} \label{0}
Suppose that $T$ is a positive closed current of bidimension $(1,1)$ on $\mathbb{P}^2$ such that $\|T\|=1$. If $\alpha \geq \frac{3}{d+5}$, $d\in \mathbb Z^+$ then there exists a curve $C$ of degree $d$ such that:
\begin{itemize}

\item $E^+_{\alpha}(T) \subset C $ 
\item or $E^+_{\alpha}(T)$ is a finite set and $|E^+_{\alpha}(T)\setminus C|\leq 2$.

\end{itemize}
\end{thm}

\section{Preliminaries}

We say that a plurisubharmonic function $u\in PSH(\mathbb C^n)$ has logarithmic growth if 
\begin{equation} \label{growth}
\gamma_u=\lim_{\|z\|\to\infty} \sup\frac{u(z)}{\log\|z\|}<\infty.  
\end{equation}

A plurisubharmonic function $u\in PSH(\mathbb C^n)$ is said to have a logarithmic pole of weight $\alpha>0$ at $p\in \mathbb C^n$ if $$u(z)=\alpha\log\|z-p\|+ O(1)$$ in an open neighborhood of $p$.

The following proposition plays an important role in our proof of the main theorem.
 
\begin{pro}[Proposition 2.1 in \cite{Com}] \label{prop1} Let $S=\{x_1,\dots,x_k \}\subset \mathbb C^n$ and $T$ be a positive closed current on $\mathbb P^n$ of bidimension $(1,1)$ with unit mass. If $u\in PSH(\mathbb C^n)$ has logarithmic growth with $\gamma_{u}$ defined as in equation \eqref{growth}, it is locally bounded outside a finite set, and $u(z)\leq\alpha_j\log\|z-x_j \|+O(1)$ for $z$ near $x_j$ where $\alpha_j>0$, $1\leq j\leq k$, then $$\sum_{j}^k\alpha_j\nu(T,x_j)\leq \gamma_u. $$
\end{pro}
 In order to construct a plurisubharmonic function $u$ as in Proposition \ref{prop1}, we will use the following result of Coman-Nivoche.
 
 \begin {thm}[Theorem 4.1 in \cite{CN}]\label{cn} Let $r$ be a positive integer and $P,Q$ be polynomials such that $S=\{z\in \mathbb C^2: P(z)=Q(z)=0 \}$ and for any $x\in S$, $ord(P,x)\geq r$,   $ord(Q,x)\geq r$  and the intersection number $(P\cdot Q)_x=r^2.$ Then $$u=\frac{1}{2r}\log(|P|^2+|Q|^2)$$ is a plurisubharmonic function in $\mathbb C^2$ with logarithmic growth, which is locally bounded outside of $S$ and has a logarithmic pole of weight $1$ at each point of $S$. 
 
 \end{thm}

\section{Proof of the main result}

%\begin{lem} Any irreducible curve of degree $d\geq 4$ may contain at most $\frac{d(d+3)}{2}$ points of $E^+_{\frac{2}{d+3}}(T).$
%\end{lem}
%\begin{proof} Suppose the contrary that an irreducible curve $C$ of degree $d$ contains $\frac{d(d+3)}{2}+1$ points of $E^+_{\frac{2}{d+3}}(T)$, say $x_1,\dots , x_{\frac{d(d+3)}{2}+1}$. By Siu's decomposition theorem $T=a[C]+R$ where $a\geq 0$ is the generic Lelong number of $T$ along $C$ and $R$ is a positive, closed  $(1,1)$ current which vanishes on $C$ generically.  By Demailly Regularization [], for any $\epsilon >0$  there exists a positive, closed  $(1,1)$ current $R'$ such that $\|R'\|=\|R\|=1-ad$, $\nu(R',x)>\nu(R,x)-\epsilon$ and $R'$ is smooth near any point at which $R$ has $0$ Lelong number. Thus by $[]$, $R'\wedge [C]$ is well-defined. Now it follows that $$ \left(\frac{d(d+3)}{2}+1\right) \left(\frac{2}{d+3}-a-\epsilon \right)<   \sum_{i=1}^{\frac{d(d+3)}{2}+1} \nu(R',x_i)\nu([C],x_i )\leq \int_{\mathbb P^2} R'\wedge [C]=(1-ad)d.$$  As $\epsilon$ is  arbitrary, $a(3d-d^2+2)\geq \frac{4}{d+3}$, which is a contradiction since left hand side is negative when $d\geq4$.
%\end{proof}

Let $d\geq 2, S=\left\{x_1,\dots, x_{\frac{d(d+3)}{2}+3} \right \}$ be a subset of $\mathbb P^2$ and $\chi$ be the subspace of $\mathbb C[X,Y,Z]$ containing all homogeneous polynomials $P(X,Y,Z) $ of degree $2d$ such that $P$ vanishes on all points of $S$ and the first derivatives of $P$ vanish on $\left\{x_1,\dots, x_{\frac{d(d+3)}{2}-3} \right \}$. Then we can estimate the dimension of $\chi$ by 
$$dim (\chi)\geq \binom{2d+2}{2}-\left( \frac{d(d+3)}{2}+3+d(d+3)-6\right)=\frac{(d-1)(d-2)}{2}+3.  $$

Note that through any set of $\displaystyle \binom {d+2}{2}-1=\frac{d(d+3)}{2}$ points in $\mathbb P^2$, there passes a (possibly reducible) curve of degree $d$.  
Let $C_1$ and $C_2$ be curves of degree $d$ passing through the points $$  x_1,\dots , x_{\frac{d(d+3)}{2}} \;\text{and}\; x_1,\dots, x_{\frac{d(d+3)}{2}-3}, x_{\frac{d(d+3)}{2}+1},x_{\frac{d(d+3)}{2}+2}, x_{\frac{d(d+3)}{2}+3}$$ respectively.  

Throughout the paper, by abuse of notation, for an algebraic curve $C\subset \mathbb P^2$ its homogenous defining polynomial will also be denoted by $C$.

\begin{lem}\label{a} Suppose that the curves $C_1$ and $C_2$ are irreducible such that $|C_1\cap C_2\cap S|\leq  \frac{d(d+3)}{2}-2.$ Then $\chi$ contains two polynomials $P,Q$ such that the common zero set of $P$ and $Q$ is discrete and contains $S$.
\end{lem}

\begin{proof} Let $P_1=C_1C_2$. Then $P_1\in \chi.$ Since $\displaystyle  dim(\chi)\geq \binom{d-1}{2}+3=:m$, we can find $P_2,\dots, P_m\in\chi$ such that $P_1,\dots, P_m$ are linearly independent. If neither $C_1$ nor $C_2$ divide $P_i$ for some $  i\in\{2,\dots, m\}$ then we can take $P=P_1$ and $Q=P_i$.  If both $C_1$ and $C_2$ divide $P_i$ for some 
$i\in\{2,\dots, m\}$, then  $P_1=\lambda P_i$ for some constant $\lambda$ which contradicts to the fact that $P_1$ and $P_i$ are linear independent.  Suppose now that $C_1$ divides $P_i$ and $C_2$ divides $P_j$ for some $i,j\in\{2,\dots, m\}$  such that $i\neq j$. Then neither $C_1$ nor $C_2$ divides $P_i+P_j$. Thus we can take $P=P_1$ and $Q=P_i+P_j\in \chi$ and we are done. The only remaining case is when one of the curves, say $C_1$,  divides all $P_2,\dots, P_m$ while the other curve,  $C_2$,  does not divide any of them. Let $P_2=C_1D_2, \dots , P_m=C_1D_m$ for some curves $D_2, \dots ,\ D_m$ of degree $d$. Note that $C_1$ may have at most $\displaystyle \binom{d-1}{2}$ singular points. Thus at least $\displaystyle \frac{d(d+3)}{2}-3-\binom{d-1}{2}=3d-4$  points of the set $\left \{x_1, \dots , x_{\frac{d(d+3)}{2}-3} \right \}$   are regular points of $C_1$. Say $x_1, \dots, x_{3d-4}$ are regular points of $C_1.$ The initial assumption $\displaystyle |C_1\cap C_2\cap S|\leq  \frac{d(d+3)}{2}-2$ implies that $C_1$ may contain at most one point of the  set   $\left\{x_{\frac{d(d+3)}{2}+1},x_{\frac{d(d+3)}{2}+2}, x_{\frac{d(d+3)}{2}+3} \right\},$ say   $x_{\frac{d(d+3)}{2}+1}$.  It follows that  $C_2\cap D_i\supset \Gamma'$ where $\Gamma'=\left \{x_1,\dots, x_{3d-4}, x_{\frac{d(d+3)}{2}+2}, x_{\frac{d(d+3)}{2}+3} \right \}$ for all $i=2,\dots \, m.$

Note that $C_2$ and $D_2$ have no common component; otherwise, we would have $C_2=D_2$, which would imply $P_1=P_2$. Let $ C_2\cap D_2=\Gamma$. By B\'ezout's Theorem  $|\Gamma|=d^2$ and $\Gamma=\Gamma' \cup \Gamma^{''} $ where $|\Gamma^{''}|=d^2-3d+2$, counting with multiplicity. By Cayley-Bacharach Theorem (Theorem CB7 in \cite{EGH}), the dimension of the vector space of curves of degree $d$ passing through the points of $\Gamma'$ modulo the ones passing through $\Gamma$ is given by $d^2-3d+2-h_{\Gamma^{''}}(d-3)$ where $h_{\Gamma^{''}}(d-3)=dim\; \mathcal O_{\mathbb P^2}(d-3)-dim\; \mathcal O_{\Gamma^{''}}(d-3)$. Let $R$ be a polynomial of degree $d-3$ vanishing on $\Gamma^{''}$. Then $|R\cap C_2|\geq d^2-3d+2>d(d-3)$. Hence by B\'ezout's Theorem, $R=C_2$ which contradicts to the fact that $R$ is of degree $d-3$. Thus $\mathcal O_{\Gamma^{''}}(d-3)$ is empty, $\displaystyle h_{\Gamma^{''}}(d-3)=dim\; \mathcal O_{\mathbb P^2}(d-3)=\binom{d-1}{2}$ and 
the dimension of the vector space of curves of degree $d$ passing through the points of $\Gamma'$ modulo the ones passing through $\Gamma$ is  $\displaystyle d^2-3d+2-\binom{d-1}{2}.$ Note that the curves $D_3, \dots, D_m$ vanish on $\Gamma'$. Since $\displaystyle m-2=\binom{d-1}{2}+1>d^2-3d+2-\binom{d-1}{2}$, a linear combination $a_3D_3+\dots +a_mD_m$ vanishes on $\Gamma$. By $AF+BG$ Theorem (Max Noether's Fundamental Theorem in \cite{Ful2}) , $ a_3D_3+\dots +a_mD_m=AC_2+BD_2$ for some $A,B.$ But this contradicts to the linear independence of $P_1, \dots , P_m.$
\end{proof}

\begin{proof}[Proof of Theorem \ref{0}] We will use induction on $d$.  For $d=1,2,3$, the result follows from \cite{Com} and \cite{KY}.  Suppose that the theorem holds for any $d'<d$. First suppose that $E_{\beta}(T)$ contains a curve $C_a$ of degree $a$ for some $ \displaystyle \beta>\frac{3}{d+5}$. We will show that $ E^+_{\frac{3}{d+5}}(T)\subset C$ for some curve $C$ of degree at most $d$.  By  Siu’s decomposition theorem, \cite{Siu}, $\displaystyle T=\beta [C_a]+R$, where $R$ has zero Lelong number generically on $C_a$. Since $\|T\|=1$ it follows that  $\displaystyle a\leq \frac{1}{\beta}<\frac{d+5}{3}$ and $\displaystyle \|R\|=\|T\|-\beta \|[C_a] \|=1-\beta a<1-\frac{3a}{d+5}$. Hence $\displaystyle E^+_{\frac{3}{d+5}}(T)=C_a\cup E^+_{\frac{3}{d+5}}(R)$ and $ \displaystyle E^+_{\frac{3}{d+5}}(R)\subset E^+_{\frac{3}{d+5-3a}}\left(\frac{R}{\|R\|}\right).$ By induction hypothesis $\displaystyle \left|E^+_{\frac{3}{d+5-3a}}\left(\frac{R}{\|R\|}\right)\setminus C'\right |\leq 2$ for some  curve $C'$ of degree $d-3a$. Hence, $E^+_{\frac{3}{d+5}}(T)$ is contained in a curve of degree at most $d-a$.

Now assume that $E_{\beta}(T)$ is $0$ dimensional for all $\beta>\alpha.$ Since $E_{\beta}(T)$ is an algebraic subvariety of $\mathbb P^2$,  $E_{\beta}(T)$ is finite.  Then $ \displaystyle  E_{\alpha}^+(T)= \bigcup_{n=1}^{\infty}E_{\alpha+\frac{1}{n}}(T)$ is at most countable. Note that $ E_{\frac{3}{d+5}}(T)\supset E_{\alpha}^+(T) $  is an algebraic subvariety of $\mathbb P^2$ of  dimension at most $1$. Hence if $E_{\alpha}^+(T)$ contains infinitely many points then $E_{\frac{3}{d+5}}(T)$ contains an algebraic curve $C_a$ of degree $a$ and result follows as above.

Finally we assume that $E_{\alpha}^+(T)$ is finite and $|E_{\alpha}^+(T)\setminus C|>2$ for any curve $C$ of degree $d$. Note that through any set of $\displaystyle \binom {d+2}{2}-1=\frac{d(d+3)}{2}$ points in $\mathbb P^2$, there passes a (possibly reducible) curve of degree $d$.  Thus we may assume that there exists a subset $S=\left\{x_1,\dots, x_{\frac{d(d+3)}{2}+3} \right \}\subset E^+_{\alpha}(T)$ such that $x_1,\dots, x_{\frac{d(d+3)}{2}} \in C_1$ and the points $x_{\frac{d(d+3)}{2}+1},  x_{\frac{d(d+3)}{2}+2}, x_{\frac{d(d+3)}{2}+3}   $ are not on $C_1$ for some curve $C_1$ of degree $d$.

 Let  $C_2$ be a curve of degree $d$ passing through the points $$  x_1,\dots, x_{\frac{d(d+3)}{2}-3}, x_{\frac{d(d+3)}{2}+1},x_{\frac{d(d+3)}{2}+2}, x_{\frac{d(d+3)}{2}+3}.$$

\textbf{Case 1:} $C_1$ and $C_2$ are irreducible.  Suppose first that  $C_2$ contains one more point of the set $S$, say $x_{\frac{d(d+3)}{2}-2}.$ Define $ u=\frac{1}{2}\log (|C_1|^2+|C_2|^2)$. By Propositon \ref{prop1}, $$\left( \frac{3}{d+5}\right)\left(\frac{d(d+3)}{2}-2\right)<\sum_{i=1}^{\frac{d(d+3)}{2}-2}\nu(T,x_i) \leq d.$$ 

This inequality fails if $d\geq 4$. Thus we may assume that $\displaystyle |C_1\cap C_2\cap S|=\frac{d(d+3)}{2}-3.$

 Let $P,Q$ be as in Lemma \ref{a}. By Theorem \ref{cn}, $\displaystyle u=\frac{1}{2}\log(|P|^2+|Q|^2)$ is a plurisubharmonic function which has logarithmic poles at the points of $S$ of total weight $d(d+3)$, locally bounded outside a finite set and $\gamma_u=2d.$ Then by Proposition \ref{prop1} we get a contradiction: 
\begin{eqnarray*} d(d+3)\frac{3}{d+5}<\sum_{j=1}^{\frac{d(d+3)}{2}+3}\nu(u,x_j)\nu(T,x_j)\leq 2d.\end{eqnarray*}

 \textbf{Case 2:} $C_1$ or $C_2$ is reducible. Without loss of generality, we may assume that $C_1$ is reducible. Say $C_1=C_aC_{a'}$ where $C_a$ is a curve of degree $a$ and $C_{a'}$ is an irreducible curve of degree $d-a$.   Let $P_1=C_1C_2\in \chi, P_2\in \chi$ be as in the proof of Lemma \ref{a}. If $P_1$ and $P_2$ do not have a common irreducible factor then $\displaystyle u=\frac{1}{2}\log (|P_1|^2+|P_2|^2)$ is a plurisubharmonic function in $\mathbb C^2$ which is locally bounded outside of a discrete set. Moreover $u$ has logarithmic pole of weight $2$ at each  points $x_1, \dots , x_{\frac{d(d+3)}{2}-3}$ and logarithmic pole of weight $1$ at each of the remaining points of $S$. It follows  from Propositon \ref{prop1} that 
 
 \begin{eqnarray*} d(d+3)\frac{3}{d+5}<\sum_{j=1}^{\frac{d(d+3)}{2}+3}\nu(u,x_j)\nu(T,x_j)\leq 2d\end{eqnarray*}
 which gives a contradiction.

 Thus we may assume that $C_{a'}, C_a, C_2$ or a component of $C_a$ or $C_2$ divides $P_2$. Suppose first that any irreducible curve of degree $s<d$ may  contain at most $\displaystyle \left\lceil{ \frac{s(d+3)}{2}}\right \rceil$ points of $S$ where $ \lceil \cdot \rceil$ denotes the ceiling function.  If   $C_{a'}$ divides $P_2$ then $P_2=C_{a'}\gamma$ where deg $\gamma=d+a.$ If $\gamma$ and $C_aC_2$ do not share a common factor then we define a plurisubharmonic function $\displaystyle u=\frac{1}{2}\log(|\gamma|^2+|C_aC_2|^2)$ which has logarithmic poles at some points of $S$ and is locally finite outside of a discrete set. By the assumption, $C_{a'}$ may contain at most $\displaystyle \left\lceil{ \frac{(d-a)(d+3)}{2}}\right\rceil $ points of $S$. Thus the total weight of $u$ at the points of $S$ is at least $ \displaystyle d(d+3)-\frac{(d-a)(d+3)}{2}-1=\frac{(d+3)(d+a)}{2}-1.$ By Proposition \ref{prop1} we get that  $$\left(\frac{(d+3)(d+a)}{2}-1\right)\frac{3}{d+5}   <\sum_{j=1}^{\frac{d(d+3)}{2}+3}\nu(u,x_j)\nu(T,x_j)\leq d+a.$$ 
 This inequality fails when $d\geq 3$ and gives a contradiction.

 If $\gamma$ and $C_aC_2$ share a common irreducible factor, say $\gamma'$, then we obtain a contradiction in a similar way using $\tilde\gamma$ and $\tilde C$ where $\gamma=\gamma'  \tilde \gamma$ and $C_aC_2=\gamma'\tilde C.$ A similar argument works if a component of $C_a$ or $C_2$ divides $P_2$.

Hence we  may  assume that there is an irreducible curve $C_s$ of degree $1 \leq s<d$ which contains at least $ \displaystyle  \left\lceil{\frac{s(d+3)}{2}}\right\rceil+1$ points of $S$, say $ x_1,\dots , x_{\left\lceil{\frac{s(d+3)}{2}}\right\rceil+1}$.  Let  $T=m[C_s]+R$ be the Siu's decomposition \cite{Siu} where $m\geq 0$ and $R$ has zero Lelong number generically on $C_s.$   By Demailly's Regularization \cite {Dem}, for any $\epsilon >0$  there exists a positive, closed  $(1,1)$ current $R'$ such that $  \|R'\|=\|R\|=1-ms$, $\nu(R',x)>\nu(R,x)-\epsilon$ and $R'$ is smooth generically on $C_s$. Thus  $ R'\wedge [C_s]$ is well-defined. It follows that 

\begin{eqnarray*}  s(1-ms)&=&\int_{\mathbb P^2}R'\wedge [C_s]\geq \sum_{i=1}^{\left\lceil{\frac{s(d+3)}{2}}\right\rceil+1}\nu(R'\wedge [C_s],x_i)\\ 
	&\geq& \sum_{i=1}^{\left\lceil{\frac{s(d+3)}{2}}\right\rceil+1}\nu(R',x_i)\nu([C_s],x_i)> \sum_{i=1}^{\left\lceil{\frac{s(d+3)}{2}}\right\rceil+1} \left(   \frac{3}{d+5} -mr_i-\epsilon  \right)r_i
	 \end{eqnarray*}
where $r_i\geq 1$ is the multiplicity of $C_s$ at $x_i$. As $\epsilon >0$ is arbitrary, 
$$m\left(  \sum_{i=1}^{\left\lceil{\frac{s(d+3)}{2}}\right\rceil+1} r_i^2-s^2 \right)\geq \frac{3}{d+5}\sum_{i=1}^{\left\lceil{\frac{s(d+3)}{2}}\right\rceil+1} r_i-s.$$

Note that $$\frac{3}{d+5}\sum_{i=1}^{\left\lceil{\frac{s(d+3)}{2}}\right\rceil+1} r_i-s\geq \frac{3}{d+5}\left(  \frac{s(d+3)}{2}+1\right)-s=\frac{s(d-1)+6}{2d+10}>0.$$
Since $m\geq 0$, it follows that $ \displaystyle \sum_{i=1}^{\left\lceil{\frac{s(d+3)}{2}}\right\rceil+1} r_i^2-s^2 >0.$ Thus 
$$\displaystyle m\geq \frac{ \frac{3}{d+5}   \displaystyle \sum_{i=1}^{\left\lceil{\frac{s(d+3)}{2}}\right\rceil+1} r_i-s }{  \displaystyle \sum_{i=1}^{\left\lceil{\frac{s(d+3)}{2}}\right\rceil+1} r_i^2-s^2}.$$

By degree-genus formula \cite{Ful2}, $$ \sum_{i=1}^{\left\lceil{\frac{s(d+3)}{2}}\right\rceil+1} r_i(r_i-1)\leq (s-1)(s-2).$$
Let $ \displaystyle P= \sum_{i=1}^{\left\lceil{\frac{s(d+3)}{2}}\right\rceil+1} r_i$ and $\displaystyle  Q= \sum_{i=1}^{\left\lceil{\frac{s(d+3)}{2}}\right\rceil+1} r_i^2.$ Then $ Q-P\leq (s-1)(s-2)=s^2-3s+2.$

It follows from this inequality that 

\begin{eqnarray*} 
m-\frac{1}{d+5}&\geq& \frac{ \frac{3}{d+5}   \displaystyle \sum_{i=1}^{\left\lceil{\frac{s(d+3)}{2}}\right\rceil+1} r_i-s }{  \displaystyle \sum_{i=1}^{\left\lceil{\frac{s(d+3)}{2}}\right\rceil+1} r_i^2-s^2}-\frac{1}{d+5}=\frac{3P-sd-5s-Q+s^2}{(d+5)(Q-s^2)}	\\ &\geq& \frac{2P-2s-2-sd}{(d+5)(Q-s^2)}\geq  \frac{s}{(d+5)(Q-s^2)}>0.
\end{eqnarray*}
Last inequality follows from the fact that $r_i\geq 1$ for all $i$.

Hence we conclude that $T=m[C_s]+R$ where $\displaystyle m> \frac{1}{d+5}$ and $ \displaystyle \|R \|=1-ms< \frac{d+5-s}{d+5}$. Thus we have that $$E^+_{\frac{3}{d+5}}(T)\subset C_s\cup E^+_{\frac{3}{d+5}}(R)\subset C_s\cup E^+_{\frac{3}{d+5-s}}\left(\frac{R}{\|R\|} \right).$$ 
By induction hypothesis there is a curve $C_{d-s}$ of degree $d-s<d$ such that  $\displaystyle \left |E^+_{\frac{3}{d+5-s}}\left(\frac{R}{\|R\|} \right)\setminus C_{d-s}\right |\leq 2.$ Hence $ \left|E^+_{\frac{3}{d+5}}(T)\setminus (C_s\cup C_{d-s})\right|\leq 2$ and we are done.

\end{proof}

\begin{rem} We note that if the points of $E^+_{\alpha}(T)$  are not overly concentrated on curves of lower degree, in the sense that
$$
m_s(E^+_{\alpha}(T))\leq \frac{s(d+3)}{2},
$$
where $m_s(E^+_{\alpha}(T))$ denotes the maximum number of points of $E^+_{\alpha}(T)$ lying on an algebraic curve of degree $s$, then the conclusion of Theorem \ref{0} for $\displaystyle \alpha=\frac{2}{d+3}$ follows readily by an argument along the same lines of the proof of Theorem \ref{0}. Hence we have the open question: Suppose that  $T$ is a positive closed  current of bidimension $(1,1)$ on $\mathbb P^2$ with unit mass and $d\geq 1$. Does there exist a degree $d$ curve $C$ such that $ |E^+_{\frac{2}{d+3}}(T)\setminus C |\leq 2 ?$ The following example shows that $\displaystyle \alpha=\frac{2}{d+3}$ would be a sharp bound if there is a positive answer to this question.
	
	\end{rem}

\begin{ex} 
Consider lines $L_1, \dots L_{d+3}$  intersecting generically in $\mathbb P^2$ and $\displaystyle T=\frac{1}{d+3}\displaystyle  \sum_{i=1}^{d+3}[L_i]$. Then the set $E_{\frac{2}{d+3}}(T)$ contains $\displaystyle \binom{d+3}{2}=\frac{(d+3)(d+2)}{2}$ many points of intersection of these lines. We claim that there does not exist a curve of degree $d\geq 1$ which contains $\displaystyle \frac{(d+3)(d+2)}{2}-2$ many points of intersections. We will use induction on $d$. For $d=1$ and $d=2$ this was proved in \cite{Com} and \cite{KY}  respectively. Suppose now the claim holds for $d-1$, $d\geq 3$ and $C$ be a  curve of degree $d$ containing  $\displaystyle\frac{(d+3)(d+2)}{2}-2$ many points of intersections.  Let $L=  \displaystyle \bigcup_{i=1}^{d+3} L_i.$  The total intersection number at the points of $L\cap C$ is at least $ d^2+5d+2>d(d+3).$ By B\'ezout's Theorem, one of the lines, say $L_1$, is an irreducible component of $C$. Let $C=C_1\cup L_1$ where $C_1$ is a curve of degree $d-1$ which contains at least $\displaystyle \frac{(d+3)(d+2)}{2}-2-(d+2)=\binom{d+2}{2}-2$ points of intersections. This contradicts to the induction hypothesis. 

\end{ex}

\section{An application: Local integrability of  potentials}
Let $\omega$ be the standard Fubini–Study form on $\mathbb P^2$ and $T$ be a positive closed $(1, 1)$-current on $\mathbb P^2$ with unit mass.  Since $T$ is cohomologous to $\omega$, $T=\omega +dd^c\varphi$ for some quasi-plurisubharmonic
(potential) function $\varphi$.  It is well known that Lelong numbers satisfy $\nu(T,a)=\nu(\varphi, a)$ at any $a\in \mathbb P^2,$ see \cite{Dem-book}. \\

\par In \cite{DK}, Demailly and Koll\'ar introduced the $\textit{integrability exponent}$ of a (quasi) plurisubharmonic function $\varphi$ at $a$ as $$ c(\varphi,a)= \sup\{ c>0:e^{-2c\varphi} \text{ is } L^1 \text{ on a neighborhood of } a\}.$$ The integrability exponent $c(\varphi,a)$ is  another way of measuring the singularity of a psh function $\varphi$ at a point $a$.  A classical theorem of Skoda \cite{Sk} states that if $\nu(\varphi,a)<2$, then $e^{-\varphi}$ is integrable in a neighborhood of $a$ with respect to the Lebesgue measure. This theorem implies that $c(\varphi,a)\geq \nu(\varphi,a)^{-1}$. Combining this lower bound with our main result Theorem \ref{0}, we obtain the following quantitative lower bound for the integrability exponent.

\begin{cor} Let $T$ be a positive closed $(1, 1)$-current on $\mathbb P^2$ with unit mass and $\varphi$ be its quasi-plurisubharmonic potential. For any $d\geq 1$, there exists a curve $C$ of degree $d$ and a set $P$ of $2$ points in $\mathbb P^2$ such that for any $x\in \mathbb P^2\setminus (C\cup P)$ $$c(\varphi,x)\geq \frac{d+5}{3}.$$ That is $e^{-c\varphi}\in L^1_{loc}(\mathbb P^2\setminus C\cup P)$ for any $c<\frac{2(d+5)}{3}.$
	\end{cor}

\noindent \textbf{Acknowledgments.} The  authors are supported by T\"UB\.ITAK 1001 Proj. No 124F144.

\end{document}